\documentclass[11pt]{article}
\usepackage{amsmath}
\usepackage{amsthm}
\usepackage{amssymb}
\usepackage{mathptmx}
\usepackage{concrete}
\def\bel#1{\begin{equation}\label #1}

\def\l{\lambda}

\def\D{\Delta}

\def\nd{\noindent}
\def\p{\partial}

\def\hf{\hfill{$\Box$}}
\def\<{\leq}
\def\>{\geq}

\begin{document}

\title{\bf Gradient estimate of a Neumann eigenfunction on a compact manifold with boundary}
\author{Jingchen Hu,\ Yiqian Shi and Bin Xu}
\date{}
\maketitle

\nd{\small {\bf Abstract.}\ \ Let $e_\l(x)$ be a Neumann
eigenfunction with  respect to the positive Laplacian $\Delta$ on
a compact Riemannian manifold $M$ with boundary such that
$\Delta\, e_\l=\l^2 e_\l$ in the interior of $M$ and the normal
derivative of $e_\l$ vanishes on the boundary of $M$. Let
$\chi_\lambda$ be the unit band spectral projection operator
associated with the Neumann Laplacian and $f$ a square integrable
function on $M$. We show the following gradient estimate for
$\chi_\lambda\,f$ as $\lambda\geq 1$: $\|\nabla\ \chi_\l\
f\|_\infty\leq C\left(\l \|\chi_\l\
 f\|_\infty+\l^{-1}\|\Delta\ \chi_\l\  f\|_\infty\right)$, where $C$ is
a positive constant depending only on $M$. As a
 corollary, we obtain the
gradient estimate of $e_\l$: for every $\l\geq 1$, there holds
$\|\nabla e_\l\|_\infty\leq C\,\l\, \|e_\l\|_\infty$.  }

\footnote{\hspace{-0.7cm} The second author is supported in part
by the National Natural Science Foundation of China (No. 10671096,
No. 10971104), the third author by the National Natural Science
Foundation of China (No. 10601053, No. 10871184). All of the
authors are supported in part by the Fundamental Research Funds
for the Central Universities.}

\nd {\small {\bf Mathematics Subject Classification (2000):}\
Primary 35P20; Secondary 35J05}

\nd {\small {\bf Key Words:}\ Neumann eigenfunction, gradient
estimate}

\section{Introduction}

 \quad  Let $(M,\ g)$ be an $n$-dimensional
compact smooth Riemannian manifold with smooth boundary $\p M$ and
$\Delta$ the positive Laplacian on $M$. In local coordinate chart
$x=(x_1,\,\cdots,x_n)$, $\Delta$ can be expressed by
\[\Delta=-\frac{1}{\sqrt{g}}\sum_{i,\,j}\, \p_{x_i}\left(g^{ij}\ \sqrt{g}\ \p_{x_j}\right),\]
where $(g^{ij})=\bigl(g^{ij}(x)\bigr)$ is the inverse of the
metric matrix $(g_{ij})=\bigl(g_{ij}(x)\bigr)=g(\p_{x_i},\
\p_{x_j})$, and
$\sqrt{g}=\sqrt{g(x)}:=\det\,\bigl(g_{ij}(x)\bigr)$. In this
paper, we always mean doing the summation from $1$ to $n$ when we
omit the variation domain of indices. Let $L^2(M)$ be the space of
square integrable functions on $M$ with respect to the Riemannian
density $dV=\sqrt{{\bf g}(x)}\,dx$. Let $e_1(x),\,e_2(x),\,\cdots$
be a complete orthonormal basis in $L^2(M)$ for Neumann
eigenfunctions of $\D$ such that $0=\l_1^2<\l_2^2\leq \l_3^2\leq
\,\cdots$ for the corresponding eigenvalues, where $e_j(x)$
($j=1,2,\dots$) are real valued smooth function on $M$ and $\l_j$
are nonnegative numbers. Also, let ${\bf e}_j$ denote the
projection of $L^2(M)$ onto the 1-dimensional space ${\bf C}e_j$.
Thus , an $L^2$ function $f$ can be written as
$f=\sum_{j=0}^\infty {\bf e}_j(f)$, where the partial sum
converges in the $L^2$ norm. Let $\l$ be a positive real number
$\geq  1$. We define the unit band spectral projection operator
(UBSPO) $\chi_\l$ as follows:
\[{\chi}_\l f:=\sum_{\l_j\in (\l,\,\l+1]}{\bf e}_j(f)\ .\]
We call that this $\chi_\l$ is associated with the Neumann
Laplacian on $M$.  The corresponding UBSPO $\chi_\l$, where we use
the same notion, can also be defined for both the Dirichlet
Laplacian on $M$ and the Laplacian on a closed Riemannian
manifold.

Grieser \cite{Gr} and Sogge \cite{S3} proved the following
$L^\infty$ estimate on $\chi_\l$ associated with the Dirichlet
Laplacian,
\begin{equation}
\label{equ:sup}
 ||{ \chi}_\l\  f||_\infty\leq C\l^{(n-1)/2} ||f||_2\,
\end{equation}
where $||f||_r$ ($1\leq r\leq \infty$) means the $L^r$ norm of the
function $f$ on $M$. In the whole of this paper $C$ denotes a
positive constant which depends only on $M$ and may take different
values at different places unless otherwise stated. The idea of
Grieser and Sogge is to use the standard wave kernel method
outside a boundary layer of width $C \lambda^{-1}$ and a maximum
principle argument inside that layer. Using their idea, Xu
\cite{X3} proved  estimate (\ref{equ:sup}) for the Neumann
Laplacian. On the other hand, Hart F. Smith \cite{Sm} proved a
sharp $L^2\to L^p$ estimate for $\chi_\l$ on a closed manifold
with Lipschitz metric. As a consequence,  (\ref{equ:sup}) holds
for both the Dirichlet or Neumann Laplacian provided $\dim\,M=2\
{\rm or}\ 3$.

By using the maximum principle argument and the estimate
(\ref{equ:sup}), Xu \cite{X2, X3} proved the following gradient
estimate on $\chi_\lambda$ for both the Dirichlet and Neumann
Laplacian,
\begin{equation}
\label{equ:2infgrad} ||\nabla\,  { \chi}_\l f||_\infty\leq
C\l^{(n+1)/2} ||f||_2\ .
\end{equation}
Here $\nabla$ is the Levi-Civita connection on $M$. In particular,
$\nabla f=\sum_j\, g^{ij}\partial f/\partial x_j$ is the gradient
vector field of a $C^1$ function $f$ in a local coordinate chart
$(x_1,\,\cdots,\,x_n)$, the square of whose length equals
$\sum_{i,j} g^{ij}(\p f/\p x_i)(\p f/\p x_j)$. One of his
motivation is to prove the H{\" o}rmander multiplier theorem on
compact manifolds with boundary. Seeger and Sogge \cite{SS}
firstly proved that theorem on closed manifolds by using the
parametrix of the wave kernel. Duong-Ouhabaz-Sikora \cite{DOS}
proved a general spectral multiplier theorem on closed manifolds
by the $L^2$ norm estimate of the kernel of spectral multipliers
and the Gaussian bounds for the corresponding heat kernel. As an
application, they gave an alternative proof to the H{\" o}rmander
multiplier theorem on closed manifolds by using the $L^\infty$
estimate \eqref{equ:sup} of $\chi_\l$ and the heat kernel.

By rescaling $\chi_\l f$ at the scale of $\l^{-1}$ both outside
and inside the boundary layer of width $C\l^{-1}$, for $\chi_\l$
associated with the Dirichlet Laplacian, the last two authors
\cite{SX2} obtained by elliptic a priori $C^{1,\,\alpha}$
estimates the following estimate slightly finer than
\eqref{equ:2infgrad},
\begin{equation}
\label{equ:grad2}
 \|\nabla\  \chi_\l\  f\|_\infty\leq C\left(\l \|\chi_\l\
 f\|_\infty+\l^{-1}\|\Delta\ \chi_\l\  f\|_\infty\right), \quad f\in L^2(M).
\end{equation}
See \cite[Remark 1.2]{SX2} for the argument that the above
estimate with the help of (\ref{equ:sup}) could imply the
Dirichlet case of estimate (\ref{equ:2infgrad}) by Xu. On the
other hand, an immediate consequence of our estimate
(\ref{equ:grad2}) is as follows: there exists a constant $C$ such
that for each Dirichlet eigenfunction $e_\l$, i.e., $\Delta\,
e_\l=\l^2\, e_\l$ in the interior of $M$ and $e_\l=0$ on the
boundary of $M$, we have $ \|\nabla e_\l\|_\infty\leq
C\,\l\,\|e_\l\|_\infty$. Furthermore, following the idea of Br{\"
u}ning \cite{Br} and Zelditch \cite[Theorem 4.1]{Z}, the last two
authors (\cite[Lemma 2.2]{SX2}) proved a basic geometry property
of nodal sets for Dirichlet eigenfunctions, i.e. {\it as $\l$
sufficiently large, every geodesic ball with radius $C/\l$ and
lying in the interior ${\rm Int}(M)$ of $M$ must contain at least
one zero point of a Dirichlet eigenfunction with eigenvalue
$\l^2$.} We call this the {\it equidistribution property} of a
non-trivial Dirichlet eigenfunction, using which we obtained a
two-sided gradient estimate for a non-trivial Dirichlet
eigenfunction $e_\l$,
\begin{equation}
\label{equ:grad} C^{-1}\, \l\,\|e_\l\|_\infty\leq \|\nabla
e_\l\|_\infty\leq C\,\l\, \|e_\l\|_\infty\quad \text{for all
$\l\geq 1$}.
\end{equation}
In the paper, we obtain in part the Neumann version of
results by the last two authors.\\

\nd {\bf Theorem 1.1.} {\it Let $f$ be a square integrable
function on the compact Riemannian manifold $(M,\,g)$ with
boundary $\p\,M$. Let $\chi_\l$ be the UBSPO associated with the
Neumann Laplacian. Then, for all $\l\geq 1$ and for all $f\in
L^2(M)$ , there holds
\begin{equation}
\label{equ:grad3}
 \|\nabla\ \chi_\l\  f\|_\infty\leq C\left(\l \|\chi_\l\
 f\|_\infty+\l^{-1}\|\Delta\ \chi_\l\  f\|_\infty\right).
\end{equation}
In particular, letting $f=e_\l(x)$ be an eigenfunction with
respect to the positive Neumann Laplacian  on $M$, i.e., $\Delta\,
e_\l=\l^2\, e_\l$ in the interior of $M$ and the normal derivative
of $e_\l$ vanishes on the boundary of $M$, we
obtain $$\|\nabla e_\l\|_\infty\leq C\,\l\,\|e_\l\|_\infty.$$}\\

\nd{\bf Remark 1.1.} We shall prove Theorem 1.1 directly via
maximum principle argument in Section 3. It is quite different
from that of the Dirichlet case (\ref{equ:grad2}) in \cite[Section
3]{SX2}, where the last two authors used the $C^{1,\,\alpha}$ a
priori estimate. Moreover, our maximum principle argument in this
paper would NOT go through for the Dirichlet case.  Heuristically
speaking, we should owe the success of the
maximum principle in the Neumann case to\\

\nd {\bf Fact 1} {\it If a $C^2$ function $g$ on the half real
line $[0,\ \infty)$ satisfies $g'(0)=0$, then the even extension
of $g$  is also $C^2$ on the real line $(-\infty,\ \infty)$.}\\

\nd Our failure of using the maximum principle argument in the
Dirichlet case is partially due to\\

 \nd {\bf Fact 2} {\it If a $C^2$ function $h$ on the half
real line $[0,\ \infty)$ satisfies $h(0)=0$, then the odd
extension of
$h$ to the real line $(-\infty,\ \infty)$ is NOT $C^2$ on $(-\infty,\ \infty)$ in general.}\\

Precisely speaking,  by Fact 1, we can reduce the gradient
estimate (\ref{equ:grad3}) near the boundary to the interior case,
which will be proved by the standard maximum argument combined
with the frequency dependent rescaling technique. However, Fact 2
prevents us from doing the similar thing for the Dirichlet case.\\

\nd {\bf Remark 1.2.} Theorem 1.1 strengthens the Neumann case of
estimate \eqref{equ:2infgrad} proved by Xu \cite{X3} in the sense
that it shows how the gradient estimate on a Neumann eigenfunction
depends on its supremum. In particular, the similar argument as
\cite[Remark 1.2]{SX2} shows that estimate \eqref{equ:grad3}
together with \eqref{equ:sup} imply the Neumann case of estimate
\eqref{equ:2infgrad} by Xu \cite{X3}. However,
\eqref{equ:2infgrad} is strong enough for Xu to prove his H{\"
o}rmander multiplier theorem associated with the Neumann Laplacian
on $M$. The authors' motivation is to prove the Neumann version of
the result of Shi-Xu \cite{SX2}.\\

\nd{\bf Remark 1.3.} We conjecture that {\it each Neumann
eigenfunction has the equidistribution property, i.e.  every
geodesic ball with radius $C/\l$ and lying in the interior ${\rm
Int}(M)$ of $M$ must contain at least one zero point of a Neumann
eigenfunction with eigenvalue $\l^2$.} If it were true, then we
could prove  the following lower bound estimate
$$\|\nabla\,e_\lambda\|_\infty\geq C\,\lambda\,\|e_\lambda\|_\infty$$
by a little modification of the argument in \cite[Section 2]{SX2}.
However, the idea of the proof for the equidistribution property
of a non-trivial Dirichlet eigenfunction in \cite[Section 2]{SX2}
did not go through for a Neumann eigenfunction,  because the
restriction of a Neumann eigenfunction to one of its nodal domains
only satisfies the mixed Dirichlet-Neumann boundary condition in
general.
\\

We conclude the introduction by explaining the organization of the
left part of this paper. We use the even extension and the maximum
principle to show Theorem 1.1 (\ref{equ:grad3}), which implies the
upper bound of $\nabla\, e_\l$. We also provide an alternative
proof of Theorem 1.1 by the same even extension and the
$C^{1,\,\alpha}$ a priori estimate.

\section{Estimate for gradient of eigenfunction}

\subsection{Outside the boundary layer}

\quad\ \ Recall the principle: {\it On a small scale comparable to
the wavelength $1/\l$, the eigenfunction $e_\l$ behaves like a
harmonic function.} It was developed in H. Donnelly  and C.
Fefferman \cite{DF1, DF2}  and was used extensively there. In this
section, for a square integrable function $f$ on $M$, letting
$\chi_\l$ be the UBSPO associated with the Neumann Laplacian, we
shall give a modification of this principle, which can be applied
to the Poisson equation
\[\D\,\chi_\l\, f=\sum_{\l_j\in (\l,\,\l+1]}\, \l_j^2\ {\bf
e}_j(f)\quad {\rm in}\quad {\rm Int}(M)\] with the Neumann
boundary condition satisfied by $\chi_\l f$ on $\p M$. Moreover,
in this subsection, we only do analysis outside the boundary layer
$L_{1/\l}=\{z\in M:\,d(z,\,\p M)\leq 1/\l\}$ of width $1/\l$.

Take a point $p$ with $d(p,\,\p M)\geq 1/\l$. We may assume that
$1/\l$ is sufficiently small such that there exists a geodesic
normal coordinate chart $(x_1,\dots, x_n)$ on the geodesic ball
$B(p,\,\frac{1}{2\l})$ in $M$. In this chart, we may identify the
ball $B(p,\,\frac{1}{2\l})$ with the $n$-dimensional Euclidean
ball ${\Bbb B}(\frac{1}{2\l})$ centered at the origin $0$, and
think of the function $\chi_\l f$ in $B(p,\,\frac{1}{2\l})$ as a
function in ${\Bbb B}(\frac{1}{2\l})$. Our aim in this subsection
is to show the inequality
\begin{equation}
\label{equ:outside2}
 |(\nabla\ \chi_\l\ f)(p)|\leq
C\left(\l\|\chi_\l\ f\|_{L^\infty\bigl({\Bbb
B}(\frac{1}{2\l})\bigr)}+\l^{-1}\|\D\ \chi_\l\
f\|_{L^\infty\bigl({\Bbb B}(\frac{1}{2\l})\bigr)}\right).
\end{equation}
For simplicity of notions, we rewrite $u=\chi_\l\, f$ and $v=\D\,
\chi_\l\, f$ in what follows. The Poisson equation satisfied by
$u$ in ${\Bbb B}(\frac{1}{2\l})$ can be written as
\[-\frac{1}{\sqrt{g}}\ \sum_{i,j}\,\p_{x_i}\left(g^{ij}\ \sqrt{g}\ \p_{x_j}
u\right)=v.\] Consider the following rescaled functions
\[u_\l(y)=u(y/\l)\quad \text{and}\quad v_\l(y)=v(y/\l)\quad
\text{in the ball ${\Bbb B}(1/2)$}.\] The above estimate which we
are after is equivalent to its rescaled version
\begin{equation}
\label{equ:scale}
 |(\nabla u_\l)(0)|\leq
C\left(\|u_\l\|_{L^\infty\bigl({\Bbb
B}(1/2)\bigr)}+\l^{-2}\|v_\l\|_{L^\infty\bigl({\Bbb
B}(1/2)\bigr)}\right).
\end{equation}
On the other hand, the rescaled version of the Poisson equation
has the expression,
\begin{equation}
\label{equ:PoissonScaled} \sum_{i,j}\,\p_{y_i}\,\left(g^{ij}_\l\
\sqrt{g_\l}\ \p_{y_j} u_\l\right)=-\l^{-2}\ \sqrt{g_\l}\
v_\l,\end{equation} where
\[g_{ij,\,\l}(y)=g_{ij}(y/\l),\quad
g^{ij}_\l(y)=g^{ij}(y/\l)\quad \text{and}\quad
\sqrt{g_\l}(y)=(\sqrt{g})(y/\l).\]

The last two authors \cite[Section 3.1]{SX2} proved
\eqref{equ:scale} by the interior $C^{1,\,\alpha}$ estimate (cf
Gilbarg-Trudinger \cite[Theorem 8.32, p. 210]{GT}) for a second
order elliptic equation of divergence type, where the $C^\alpha$
norm of coefficients $g^{ij}_\l\,\sqrt{g_\l}$ in the equation are
involved. In the following paragraph, we shall give a different
and more elementary proof of \eqref{equ:scale}, where we use the
maximum principle, however, the $C^{0,\,1}$ norm of coefficients
$g^{ij}_\l\,\sqrt{g_\l}$ are involved. Note that the $C^{0,\, 1}$
norms of $g^{ij}_\l\,\sqrt{g_\l}$ are uniformly bounded for all
$\l\geq 1$.

For simplicity of notions, we set
\[u_\l=\phi,\quad
h=-\l^{-2}\ \sqrt{g_\l}\ v_\l,
\quad a_{ij}=g^{ij}_\l\
\sqrt{g_\l},\quad b_i=\sum_{j=1}^n\, \frac{a_{ij}}{\p\,y_j}.\]
Then the rescaled Poission equation \eqref{equ:PoissonScaled} can
be written as
\begin{equation*}
L\phi:=\sum_{i,j}\, a_{ij}\ \frac{\p^2\,\phi}{\p y_i\p
y_j}+\sum_i\,b_i\ \frac{\p \phi}{\p y_i}=h,\quad y\in {\Bbb
B}(1/2).
\end{equation*}
We learned this maximum principle argument for proving
\eqref{equ:scale} from Brandt \cite[p. 95-6]{Bra}. Moreover, we
find that the constant-coefficient-assumption there could be
removed. The idea is to construction a new function $\phi_1$ from
$\phi$ of $n+1$ variables and apply the maximum principle to
$\phi_1$. The details are as follows. Define
\[\phi_1(y_1,\,\cdots,\,y_n;\,z_1):=
\frac{1}{2}\ \Bigl(\phi(y_1+z_1,\,y_2,\cdots,y_n)-\phi(y_1-z_1,\,y_2,\cdots,y_n)\Bigr)
\]
in the $(n+1)$-dimensional domain
\[{\cal R}=\{(y_1,\,\cdots,\,y_n;\, z_1):\,|y|<1/4,\ 0<z_1<1/4\}.\]
Writing
\[L_1=L-\mu\,\frac{\p^2}{\p y_1^2}+\mu\, \frac{\p^2}{\p
z_1^2}\quad (\mu>0),\] we observe that, for sufficiently small
$\mu$, this new operator is elliptic in the $n+1$ variables, and
satisfies
\[|L_1\ \phi_1|=|L\ \phi_1|\leq \|h\|\quad \text{in} \quad{\cal R},\]
where we denote by $\|\cdot\|$ the $L^\infty$ norm in ${\Bbb
B}(1/2)$. Choose a constant $C$ sufficiently large and depending
on the $L^\infty$ norm of coefficients $a_{ij}$ and $b_i$ so that
\[L(|y|^2)\leq 2\mu\,C\]
and introduce the comparison function
\[\overline{\phi_1}:=\frac{1}{2\mu}\,\|h\|\,\left(\frac{1}{4}\,z_1-z_1^2\right)
+16\,\|\phi\|\,\left\{|y|^2+z_1^2+C\left(\frac{1}{4}z_1-z_1^2\right)\right\}
.\]  Then we have
\begin{eqnarray*}
L_1\,
\overline{\phi_1}&=&-\|h\|+16\,\|\phi\|\,\bigl(L(|y|^2)-2\,\mu\,C\bigr)\\
&\leq&-\|h\|\leq -|L_1\,\phi_1| \quad {\rm in}\quad {\cal R}
\end{eqnarray*}
and
\[\overline{\phi_1}\geq |\phi_1|\quad \text{on the boundary}\quad
\p {\cal R}.\] Thus, by the weak maximum principle (cf
Gilbarg-Trudinger \cite[Theorem 3.1, p. 32]{GT}), we obtain
$|\phi_1|\leq \overline{\phi_1}$. This implies that
\begin{eqnarray*}
\frac{1}{2}\,\left|\phi(z_1,\,0,\cdots,0)-\phi(-z_1,\,0,\cdots,0)\right|
&\leq&  \overline{\phi_1}(0,\cdots,0,\,z_1)\\
&\leq&\frac{1}{2\mu}\,\frac{z_1}{4}\,\|h\|+16\,\|\phi\|\left(\frac{Cz_1}{4}+z_1^2\right).
\end{eqnarray*}
Dividing through by $z_1$ and letting $z_1\to 0$ yields the
desired estimate
\begin{equation}
\label{equ:outsidemax}
\left|\frac{\p\phi}{\p y_1}(0)\right|\leq
\frac{1}{8\mu}\|h\|+4C\,\|\phi\|.
\end{equation} Therefore, we complete the proof of
\eqref{equ:outside2}.

We remark that \eqref{equ:outside2} can also be proved directly by
the above maximum principle argument without doing the re-scaling.
Here we prefer to do the rescaling before proceeding to the
maximum principle argument because of the following two reasons:

\nd $\bullet$ Re-scaling makes the dependence relation of the
desired estimate on the eigenvalue $\l^2$ clear and reduce the
question to the case of a fixed scale.

\nd $\bullet$ It is convenient for reader to compare the  maximum
principle argument here with the proof via the elliptic a priori
estimate in Shi-Xu \cite{SX, SX2}.

\subsection{Inside the boundary layer}

\quad\ \ Using the notions in subsection 2.1, we are going to
prove the following estimate:
\begin{equation}
\label{equ:inside} |\nabla u (p_0)|\leq C\left(\l
\|u\|_\infty+\l^{-1}\|v\|_\infty\right) \quad \text{for all}\
p_0\in L_{1/\l},
\end{equation}
with which combining (\ref{equ:outside2}) completes the proof of
Theorem 1.1.

Since the boundary $\p M$ is a compact sub-manifold in $M$ of
codimension 1,  we can take $\l$ sufficiently large such that
there exists the boundary normal coordinate chart  on the boundary
layer $L_{3/\l}=\{p\in M:\,d(p,\,\p M)\leq 3/\l\}$ with respect to
the boundary $\p M$ (cf H{\" o}rmander \cite[p. 51]{H2}). In
particular, we have the map
\[{\cal B}: \p M\times [0,\,3/\l]\to L_{3/\l},\quad (p',\,\delta)\mapsto {\cal B}(p',\,\delta)\]
such that $\delta\mapsto {\cal B}(p',\,\delta)$, $\delta\in
[0,\,3/\l]$, is the geodesic with arc-length parameter normal to
$\p N$ at $p'$. Moreover, for each point $(p',\,\delta)\in
L_{3/\l}$, we have $0\leq \delta\leq 3/\l$ and
\[d\bigl((p',\,\delta),\,\p M\bigr)=\delta.\] Denote by ${\cal
R}(r)$ the following $n$-dimensional rectangle in ${\bf R}^n$
sitting at the origin and having size $r$,
\[{\cal R}(r)=\left\{x=(x',\,x_n)=\bigl((x_1,\dots,x_{n-1}),\,x_n\big)\in {\bf R}^n:|(x',\,0)|<r,\,
0\leq x_n\leq r\right\}.\] For a point $q$ on $\p M$, denote by
${\rm Exp}_q$ at $q$ the exponential map on the sub-manifold $\p
M$. Since $\p M$ is compact and $\l$ is sufficiently large, we may
assume the existence of the geodesic normal chart for each metric
ball of radius $3/\l$ on $\p M$.

We choose and fix a point $p_0$ in $L_{1/\l}$, and write
$p_0={\cal B}(q_0,\,\delta_0)$, where $q_0\in \p M$ and
$\delta_0\in [0,\,1/\l]$. We denote by $R(q_0,\,3/\l)$ the
rectangle in $M$ sitting at $q_0$ and having size $3/\l$,
\[R(q_0,\,3/\l)=\left\{\big({\rm Exp}_{q_0}(x'),\,x_n)\big):\,(x',\,x_n)\in {\cal R}(3/\l)\right\}.\]
In this way, we identity the rectangle $R(q_0,\,3/\l)$ in $M$
sitting at $q_0$ and containing $p_0$ with the rectangle ${\cal
R}(3/\l)$ in ${\bf R}^n$. Thus we could look at $u$ and $v$ as
functions in ${\cal R}(3/\l)$.

We recall that $u_\l$ and $v_\l$ are the corresponding rescaled
functions of $u$ and $v$, respectively,  i.e.,
\[u_\l(y)=u(y/\l),\quad v_\l(y)=v(y/\l)\quad\text{for all $y$ in ${\cal
R}(3)$}.\] To prove \eqref{equ:inside}, we need only show the
following estimate,
\begin{equation}
\label{equ:insidescaled} |(\nabla\, u)(p_0)|=|(\nabla\, u)(0,\
\delta_0)|\leq C\left(\l\, \|u\|_{L^\infty\bigl({\cal
R}(3/\l)\bigr)}+\l^{-1}\|v\|_{L^\infty\bigl({\cal
R}(3/\l)\bigr)}\right),
\end{equation}
which can be reduced to the equivalent rescaled version,
\begin{equation}
\label{equ:insidescaled} |(\nabla u_\l)(0,\ \l\,\delta_0)|\leq
C\left( \|u_\l\|_{L^\infty\bigl({\cal
R}(3)\bigr)}+\l^{-2}\|v_\l\|_{L^\infty\bigl({\cal
R}(3)\bigr)}\right),\quad 0\leq \l\,\delta_0\leq 1.
\end{equation}
where $u_\l$ and $v_\l$ are the the rescaling function of $u$ an
$v$, respectively. Observe that $u_\l$ is the solution of the
Poisson equation
\begin{equation}
\label{equ:PoissonScaled2} \sum_{i,j}\,\p_{y_i}\,\left(g^{ij}_\l\
\sqrt{g_\l}\ \p_{y_j} u_\l\right)=-\l^{-2}\ \sqrt{g_\l}\ v_\l\
\quad \text{in the interior of rectangle}\ {\cal R}(3)
\end{equation} and
satisfies the Neumann boundary condition, i.e.,
\[\frac{\p u_\l}{\p y_n}=0\ \text{
on the portion}\  \{x\in {\cal R}(3):y_n=0\}\  \text{of the
boundary}\  \p {\cal R}(3).\] We shall give two different proofs
for \eqref{equ:insidescaled}.\\

\nd {\bf  1st proof}\quad  The idea is to reduce, by Fact 1 in the
introduction and the even extension, the question to the interior
gradient estimate \eqref{equ:scale}, which has been proved by the
maximum principle in the former subsection. By the geometric
property of geodesic normal coordinate chart with respect to the
boundary $\p M$, we have
\begin{equation*}
g^{nn}(x',\ x_n)=1\  \text{and}\  g^{jn}(x',\ x_n)=0\ \text{for}\
j\not=n\  \text{in}\ R(q_0,\,3/\l),
\end{equation*}
which implies that
\begin{equation*}
g_{\l}^{nn}(y',\ y_n)=1\  \text{and}\  g_{\l}^{jn}(y',\ y_n)=0\
\text{for}\ j\not=n\  \text{in}\ {\cal R}(3).
\end{equation*}

Setting
\[
a_{ij}:=g^{ij}_\l\,\sqrt{g_\l}\quad \text{for}\quad 1\leq
i,\,j\leq n-1,\quad \ a_{nn}:=\sqrt{g_\l}\] and
\[b_i:=\sum_{j=1}^{n-1}\,\frac{\p
a_{ij}}{\p y_j}\quad \text{for}\quad 1\leq i\leq n-1,\quad
b_n=\frac{\sqrt{g_\l}}{\p y_n},
\]
we can express the Poisson equation \eqref{equ:PoissonScaled2} as
\begin{equation}
\label{equ:insidemax} \sum_{1\leq i,\,j\leq n-1}\,
a_{ij}\,\frac{\p^2\phi}{\p y_i\p y_j}+a_{nn}\,\frac{\p^2\phi}{\p
y_n^2}+\sum_{1\leq k\leq n}\,b_k\,\frac{\p \phi}{\p y_k}=h\quad
\text{in}\quad {\rm Int}\bigl({\cal R}(3)\bigr),
\end{equation}
 where $\phi=u_\l$ and $h=-\l^{-2}\ \sqrt{g_\l}\ v_\l$.

Set
\[{\cal S}(r):=\left\{y=(y',\,y_n)=\bigl((y_1,\dots,y_{n-1}),\,y_n\big)\in {\bf R}^n:|(y',\,0)|<r,\,
|y_n|\leq r\right\},\] which is the union of rectangle ${\cal
R}(r)$ and its reflection with respect to the hyperplane
$\{y_n=0\}$. We denote by ${\tilde \phi}$ the even extension onto
${\cal S}(3)$ of the function $\phi$ defined on ${\cal R}(3)$,
i.e.,
\[{\tilde \phi}(y',\ y_n)=
\begin{cases}
\phi(y',\ y_n) & \text{if}\quad y_n\geq 0\\
\phi(y',\ -y_n) & \text{if}\quad y_n< 0.
\end{cases}
\]
We do the even extension to $h$ and the coefficients
$a_{ij},\,a_{nn},\,b_i$ for $1\leq i,\,j\leq n-1$, and denote the
corresponding extension functions on ${\cal S}(3)$ by
\[{\tilde h},\quad {\widetilde {a_{ij}}},\quad  {\widetilde
{a_{nn}}},\quad {\widetilde {b_i}}.\] However we do the odd
extension to $b_n$,
\[{\widetilde {b_n}}(y',\ y_n)=
\begin{cases}
b_n(y',\ y_n) & \text{if}\quad y_n\geq 0\\
-b_n(y',\ -y_n) & \text{if}\quad y_n< 0.
\end{cases}
\]
We shall see soon that the possible discontinuity of ${\widetilde
{b_n}}$ on the portion ${\cal S}(3)\cap \{y_n=0\}$ would not cause
any trouble.

Thus, we obtain the following Poisson equation about ${\tilde
\phi}$ with continuous coefficients
\[\sum_{1\leq i,\,j\leq n-1}\,
{\widetilde {a_{ij}}}\,\frac{\p^2{\tilde \phi}}{\p y_i\p
y_j}+{\widetilde {a_{nn}}}\,\frac{\p^2{\tilde \phi}}{\p
y_n^2}+\sum_{1\leq k\leq n}\,{\widetilde {b_k}}\,\frac{\p {\tilde
\phi}}{\p y_k}={\tilde h}\quad \text{in}\quad {\rm Int}\bigl({\cal
S}(3)\bigr)
\]
except that $\widetilde{b_n}$ is bounded and possibly
discontinuous on the portion ${\cal S}(3)\cap \{y_n=0\}$. By Fact
1 in the introduction, which can be proved by simple calculus
computation, we know that ${\tilde \phi}$ is $C^2$ in ${\cal
S}(3)$. The only point which we should take care of is whether
${\widetilde {b_n}}\,\frac{\p {\tilde \phi}}{\p y_n}$ is an even
continuous function in ${\cal S}(3)$ with respect to $y_n$.
However, by the extension ${\tilde \phi}$ of $\phi$ and the
Neumann boundary condition, i.e., $\frac{\p \phi}{\p y_n}=0$ on
$\{y_n=0\}\cap {\cal S}(3)$, $\frac{\p {\tilde \phi}}{\p y_n}$ is
an odd $C^1$ function vanishing on the portion ${\cal S}(3)\cap
\{y_n=0\}$. Since ${\widetilde {b_n}}$ is a bounded function in
${\cal S}(3)$ and is odd with respect to $y_n$,  ${\widetilde
{b_n}}\,\frac{\p {\tilde \phi}}{\p y_n}$ is a continuous function
being even with respect to $y_n$ in ${\cal S}(3)$. Moreover,
${\widetilde {b_n}}\,\frac{\p {\tilde \phi}}{\p y_n}$ vanishes on
the portion $\{y_n=0\}\cap {\cal S}(3)$. Therefore, we have
reduced the proof of \eqref{equ:insidescaled} to the estimate for
${\tilde \phi}$ at the interior point $(0,\,\delta_0\,\l)$ of
${\cal S}(3)$ similar to \eqref{equ:outsidemax} in the former
subsection. The fact that the bounded coefficient $b_n$ is
possibly discontinuous on the portion $\{y_n=0\}\cap {\cal S}(3)$
does not bring us any trouble of applying the weak maximum
principle. See Gilbarg-Trudinger \cite[(3.3), p.31]{GT} and the
related comments.  \hf\\

\nd {\bf 2nd proof}\quad The idea is to use the same even
extension as above and the interior $C^{1,\alpha}$ estimate
Gilbarg-Trudinger \cite[Theorem 8.32, p. 210]{GT}. Denote by
${\tilde g}$ the even extension of the Riemannian metric $g$ on
${\cal R}(3/\l)$ onto ${\cal S}(3/\l)$. Then ${\tilde g}$ is a
Lipschitz metric on ${\cal S}(3/\l)$ with $C^{0,\,1}$ norm bounded
by the $C^1$ norm of $g$. Denote the even extension of $u$ and $v$
on ${\cal S}(3/\l)$ by ${\tilde u}$ and ${\tilde v}$,
respectively. We claim that ${\tilde u}$ is a weak solution of the
following Poisson equation
\[-\frac{1}{\sqrt{\tilde g}}\ \sum_{i,j}\,\p_{x_i}\left({\tilde g}^{ij}\ \sqrt{\tilde g}\ \p_{x_j}
{\tilde u}\right)={\tilde v}\quad \text{in}\quad {\rm
Int}\bigl({\cal S}(3/\l)\bigr).
\]
That is, for each smooth function $\psi$ compactly supported in
${\rm Int}\bigl({\cal S}(3/\l)\bigr)$, the following integral
equality holds
\begin{equation}
\label{equ:weaksol} \int_{{\rm Int}\bigl({\cal S}(3/\l)\bigr)}\,
\sum_{i,j}\, {\tilde g}^{ij}\  \p_{x_i} {\tilde u}\ \p_{x_j}\psi\
dx=\int_{{\rm Int}\bigl({\cal S}(3/\l)\bigr)}\,(-{\tilde v})\
\psi\ dx.
\end{equation}
Actually, since $\p {\tilde u}/\p x_n=0$ on ${\cal
S}(3/\l)\cap\{x_n=0\}$, we find by the Green formula on Riemannian
manifolds and the even extension of $u$ and $v$,
\begin{eqnarray*}
\int_{{\rm Int}\bigl({\cal S}(3/\l)\bigr)\cap \{x_n>0\}}\,
\sum_{i,j}\, {\tilde g}^{ij}\  \p_{x_i} {\tilde u}\ \p_{x_j}\psi\
dx&=&\int_{{\rm Int}\bigl({\cal S}(3/\l)\bigr)\cap
\{x_n>0\}}\,\Delta\,{\tilde u}\ \psi\ dx\\
&=&\int_{{\rm Int}\bigl({\cal S}(3/\l)\bigr)\cap
\{x_n>0\}}\,\Delta\,u\ \psi\
dx\\
&=& \int_{{\rm Int}\bigl({\cal S}(3/\l)\bigr)\cap
\{x_n>0\}}\,(-v)\ \psi\ dx\\
&=&\int_{{\rm Int}\bigl({\cal S}(3/\l)\bigr)\cap
\{x_n>0\}}\,(-{\tilde v})\ \psi\ dx.
\end{eqnarray*}
Using the change of variable $x_n\mapsto -x_n$ and  the above
equality, we obtain
\[
\int_{{\rm Int}\bigl({\cal S}(3/\l)\bigr)\cap\{x_n<0\}}\,
\sum_{i,j}\, {\tilde g}^{ij}\  \p_{x_i} {\tilde u}\ \p_{x_j}\psi\
dx=\int_{{\rm Int}\bigl({\cal S}(3/\l)\bigr)\cap\{x_n<0\}}\,
(-{\tilde v})\ \psi\ dx,
\]
where we also use $g^{in}=0$ for all $i\not=n$.  Summing these two
equality yields \eqref{equ:weaksol}. Recall that coefficients
${\tilde g}^{ij}$ is Lipschitz and ${\tilde v}$ is continuous on
${\cal S}(3/\l)$. On the other hand, we have the rescaled version
of equation \eqref{equ:weaksol}, i.e., for each smooth function
$\psi$ compactly supported in ${\rm Int}\bigl({\cal S}(3)\bigr)$
\[
\int_{{\rm Int}\bigl({\cal S}(3)\bigr)}\, \sum_{i,j}\, {\tilde
g}^{ij}_\l\  \p_{y_i} {\tilde u}_\l\ \p_{y_j}\psi\ dx=\int_{{\rm
Int}\bigl({\cal S}(3)\bigr)}\,(-{\tilde v}_\l)\ \psi\ dx.
\]
Thus, applying to it the interior $C^{1,\alpha}$ estimate in
Gilbarg-Trudinger \cite[Theorem 8.32, p. 210]{GT},
 we obtain that for every $0<\alpha<1$,
\[\|{\tilde u}_\l\|_{C^{1,\,\alpha}({\cal S}(2))}\leq
C\left(\|{\tilde u}_\l\|_{C^0({\cal S}(3))}+\l^{-2}\ \|{\tilde
v}_\l\|_{C^0({\cal S}(3))} \right),\] which implied the desire
estimate \eqref{equ:insidescaled}.\hf \\

\nd  {\bf Acknowledgements}\quad The last author would like to
thank Professor Qing Han, Professor Xinan Ma, Professor
Christopher D. Sogge and Professor Meijun Zhu for valuable
conversation during the course of this work. Y.S. is supported in
part by the National Natural Science Foundation of China (No.
10971104), B.X. by the National Natural Science Foundation of
China (Grant No. 11271343) and Anhui Provincial Natural Science
Foundation (Grant No. 1208085MA01). The last two authors are
supported in part by the Fundamental Research Funds for the
Central Universities (Grant No. WK0010000020 and No.
WK0010000023).

\nd {\bf Author information}\\

\nd Jingchen Hu, Wu Wen-Tsun Key Laboratory of Mathematics, USTC,
Chinese Academy of Sciences. School of Mathematical Sciences,
University of Science and Technology of China, Hefei 230026
China.\\

\nd Yiqian Shi,  Wu Wen-Tsun Key Laboratory of Mathematics, USTC,
Chinese Academy of Sciences. School of Mathematical Sciences,
University of Science and Technology of China, Hefei 230026 China.\\

\nd Bin Xu, Wu Wen-Tsun Key Laboratory of Mathematics, USTC,
Chinese Academy of Sciences. School of Mathematical Sciences,
University of Science and Technology of China, Hefei 230026 China.

\nd E-mail: bxu@ustc.edu.cn

\end{document}